\documentclass[12pt]{amsart}
\usepackage{amssymb}
\newtheorem{theorem}{Theorem}
\newtheorem{lemma}{Lemma}
\newtheorem{proposition}{Proposition}
\newtheorem{corollary}{Corollary}
\newtheorem{definition}{Definition}
\newcommand{\ffrac}[2]{{\mbox{\large$\frac{#1}{#2}$}}}

\begin{document}
\title{Higher Symmetries of the Square of the Laplacian}
\author{Michael Eastwood}
\address{School of Mathematical Sciences, University of Adelaide,
SA 5005, Australia}
\email{meastwoo@member.ams.org}
\author{Thomas Leistner}
\address{School of Mathematical Sciences, University of Adelaide,
SA 5005, Australia}
\email{tleistne@maths.adelaide.edu.au}
\thanks{This work was undertaken in preparation for and during the 2006 Summer
Program at the Institute for Mathematics and its Applications at the University
of Minnesota. The authors would like to thank the IMA for hospitality during
this time. The authors are supported by the Australian Research Council.}
\begin{abstract}
The symmetry operators for the Laplacian in flat space were recently described
and here we consider the same question for the square of the Laplacian. Again,
there is a close connection with conformal geometry. There are three main steps
in our construction. The first is to show that the symbol of a symmetry is
constrained by an overdetermined partial differential equation. The second is
to show existence of symmetries with specified symbol (using a simple version
of the AdS/CFT correspondence). The third is to compute the composition of two
first order symmetry operators and hence determine the structure of the
symmetry algebra. There are some interesting differences as compared to the
corresponding results for the Laplacian.
\end{abstract}
\dedicatory{In memory of Thomas Branson} 
\renewcommand{\subjclassname}{\textup{2000} Mathematics Subject Classification}
\maketitle

\section{Introduction} 
The second order symmetry operators for the Laplacian on~${\mathbb{R}}^n$ were
determined by Boyer, Kalnins, and Miller~\cite{boyer-kalnins-miller76}. The
higher order symmetries were found in \cite{eastwood05} and the structure of
the resulting algebra was also described. Here were prove the corresponding
results for the square of the Laplacian. The other aspect
of~\cite{boyer-kalnins-miller76}, namely the relation between second order
symmetries and separation of variables, is unclear for the square of the
Laplacian. 

We are grateful to Ernie Kalnins who suggested the square of the Laplacian as a
candidate for having interesting symmetries. We would also like to acknowledge
pertinent comments from Petr Somberg, Vladim\'{\i}r Sou\v{c}ek, and Misha
Vasiliev.

\section{Definitions and statements of results}\label{DefineandState}
We shall always work on $n$-dimensional Euclidean space ${\mathbb{R}}^n$ for
$n\geq 3$ and adopt the usual convention of writing vectors and tensors adorned
with indices, which we shall raise and lower with the standard (flat)
metric~$g_{ab}$. Let us also write $\nabla_a=\partial/\partial x^a$ for
differentiation in co\"ordinates. Then $\nabla^a=g^{ab}\nabla_b$ and the
Laplacian is given by $\Delta=\nabla^a\nabla_a$. All functions and tensors in
this article will be smooth. All differential operators will be linear with
smooth coefficients.
\begin{definition}
A differential operator ${\mathcal{D}}$ is a symmetry of $\Delta^2$ if and only
if there is another differential operator $\delta$ such that
$\Delta^2{\mathcal{D}}=\delta\Delta^2$.
\end{definition}
\noindent Obviously, any operator of the form ${\mathcal{P}}\Delta^2$ is a
symmetry of $\Delta^2$ because we can take $\delta=\Delta^2{\mathcal{P}}$.
Therefore one introduces the following equivalence relation.
\begin{definition}\label{equivalence}
Two symmetries ${\mathcal{D}}_1$ and ${\mathcal{D}}_2$ of $\Delta^2$ are
equivalent, ${\mathcal{D}}_1\sim{\mathcal{D}}_2$, if and only if
${\mathcal{D}}_1-{\mathcal{D}}_2={\mathcal{P}}\Delta^2$ for some differential
operator~${\mathcal{P}}$.
\end{definition}
\noindent Of course, this equivalence relation only effects symmetries of order
$s\geq4$. The composition of two symmetries is again a symmetry. Also,
composition is compatible with the equivalence relation, 
i.e.\ if ${\mathcal{D}}_1\sim{\mathcal{D}}_2$ and 
${\mathcal{D}}_3\sim{\mathcal{D}}_4$, then 
${\mathcal{D}}_1{\mathcal{D}}_3\sim{\mathcal{D}}_2{\mathcal{D}}_4$. 
This allows us to define an algebra:--
\begin{definition}
The algebra ${\mathcal{B}}_n$ consists of all symmetries of $\Delta^2$ 
on ${\mathbb{R}}^n$ considered modulo equivalence and with algebra operation
induced by composition.
\end{definition}
\noindent In the following we shall study this algebra and describe its
structure. To this end we need the notion of conformal Killing tensors and
their generalisations as studied in~\cite{nikitin91,nikitin-prilipko90}. 
We shall write $\phi^{(ab\cdots c)}$ for the symmetric part of a tensor
$\phi^{ab\cdots c}$. 
\begin{definition}
A conformal Killing tensor $V^{bcd\cdots e}$ is a symmetric trace-free tensor
such that 
\begin{equation}\label{confKill}
\mbox{\rm the trace-free part of }\nabla^{(a}V^{bcd\cdots e)}=0,
\end{equation}
equivalently that
$$\nabla^{(a}V^{bcd\cdots e)}=g^{(ab}\phi^{cd\cdots e)}$$
for some tensor $\phi^{cd\cdots e}$. A conformal Killing tensor with one index 
is called a conformal Killing vector. A conformal Killing tensor with no
indices is simply a constant. 
\end{definition}
\begin{definition}
A generalised conformal Killing tensor $W^{d\cdots e}$ of order 3 is a
symmetric trace-free tensor such that
$$\mbox{\rm the trace-free part of }
\nabla^{(a}\nabla^b\nabla^cW^{d\cdots e)}=0,$$
equivalently that
$$\nabla^{(a}\nabla^b\nabla^cW^{d\cdots e)}=g^{(ab}\phi^{cd\cdots e)}$$
for some tensor $\phi^{cd\cdots e}$.
\end{definition}
\noindent Though it is clear how to define a generalised conformal Killing
tensor of any order, we shall only need order $3$. This should be taken
as read for the rest of this article. 

Our main theorems on the existence and uniqueness of symmetries are as follows.
\begin{theorem}\label{Utheo}
Any zeroth order symmetry of $\Delta^2$ is of the form 
$$f\longmapsto Vf \quad\mbox{\rm for $V$ constant}.$$
Any first order symmetry of $\Delta^2$ is of the form 
$$V^b\nabla_b+\mbox{\rm lower order terms},$$
where $V^b$ is a conformal Killing vector. Any higher symmetry, say of
degree~$s$, of $\Delta^2$ is canonically equivalent to one of the form
$$V^{bcd\dots e}\nabla_b\nabla_c\nabla_d\cdots\nabla_e+
W^{d\cdots e}\Delta\nabla_d\cdots\nabla_e+\mbox{\rm lower order terms},$$
where $V^{bcd\cdots e}$ is a conformal Killing tensor of valency $s$ and
$W^{d\cdots e}$ is a generalised conformal Killing tensor of valency~$s-2$.
\end{theorem}
\begin{theorem}\label{Etheo}
Suppose that $V^{bcd\cdots e}$ is a conformal Killing tensor on
${\mathbb{R}}^n$.
Then there is a canonically defined differential operator 
$${\mathcal{D}}_V=V^{bcd\dots e}\nabla_b\nabla_c\nabla_d\cdots\nabla_e
+\mbox{\rm lower order terms}$$
that is a symmetry of $\Delta^2$. Suppose that $W^{d\cdots e}$ is a generalised
conformal Killing tensor. 
Then there is a canonically defined differential operator 
$${\mathcal{D}}_W=W^{d\dots e}\Delta\nabla_d\cdots\nabla_e
+\mbox{\rm lower order terms}$$
that is a symmetry of $\Delta^2$. 
\end{theorem}
The proof of this theorem will be given in Section \ref{Etheoprf} by using the
ambient metric construction. Here we only want to give the first and second
order symmetries. As a special case of~(\ref{DVf}), any first order symmetry is
given by
$${\mathcal{D}}f=V^a\nabla_af+\ffrac{n-4}{2n}(\nabla_aV^a)f+cf,$$
for a conformal Killing vector $V^a$ and arbitrary constant~$c$. The canonical
ones ${\mathcal{D}}_V$ are those with $c=0$. The canonical second order
symmetries are
\begin{equation}\label{DVtwo}{\mathcal{D}}_Vf=V^{ab}\nabla_a\nabla_bf
+\ffrac{n-2}{n+2}(\nabla_aV^{ab})\nabla_bf
+\ffrac{(n-2)(n-4)}{4(n+1)(n+2)}(\nabla_a\nabla_bV^{ab})f,\end{equation} for
$V^{ab}$ a conformal Killing tensor and
\begin{equation}\label{DWtwo}{\mathcal{D}}_Wf=W\Delta f-(\nabla^a W)\nabla_af
-\ffrac{n-4}{2(n+2)}(\Delta W)f,\end{equation} for $W$ a generalised conformal
Killing scalar, i.e.\ $\nabla^a\nabla^b\nabla^cW=g^{(ab}\phi^{c)}$. Of course,
there is no freedom in equivalence until we consider fourth order operators.
Hence, we can use Theorems~\ref{Utheo} and~\ref{Etheo} to write down all second
order symmetries as follows. Suppose that ${\mathcal{D}}$ is a second order
symmetry operator. According to Theorem~\ref{Utheo}, it has the form
$${\mathcal{D}}=V^{ab}\nabla_a\nabla_b+W\Delta+\mbox{lower order terms}$$
where $V^{ab}$ is a conformal Killing tensor and $W$ is a generalised conformal
Killing scalar. According to Theorem~\ref{Etheo}, however, there are
canonically defined symmetry operators of the same form, which we can subtract
to obtain a first order symmetry. Iterating this procedure we conclude that 
$${\mathcal{D}}={\mathcal{D}}_{V_2}+{\mathcal{D}}_W+{\mathcal{D}}_{V_1}
+{\mathcal{D}}_{V_0},$$
where ${\mathcal{D}}_{V_s}$ are the canonically defined differential operators
associated to conformal Killing tensors of valency $s$ and ${\mathcal{D}}_W$ is
the operator associated to a generalised conformal Killing scalar. As above, we
have explicit formulae for these operators and, of course,
${\mathcal{D}}_{V_0}f=V_0f$ for constant~$V_0$. We shall soon see that the
space of (generalised) conformal Killing tensors is finite-dimensional. (In
particular, the space of second order symmetries of $\Delta^2$ has dimension
$(n+1)(n+2)(n^2+5n+12)/12$.)

More generally, let ${\mathcal{K}}_{n,s}$ denote the vector space of conformal
Killing tensors on ${\mathbb{R}}^n$ with $s$ indices and, for $s\geq 2$, let
${\mathcal{L}}_{n,s}$ denote the vector space of generalised conformal Killing
tensors on ${\mathbb{R}}^n$ with $s-2$ indices. Reasoning as we just did for
second order symmetries, but now taking into account the equivalence
necessitated by Theorem~\ref{Utheo} in general, we conclude that any symmetry
of $\Delta^2$ may be canonically thrown into an equivalent one of the form 
$${\mathcal{D}}_{V_s}+{\mathcal{D}}_{W_s}+\cdots
+{\mathcal{D}}_{V_2}+{\mathcal{D}}_{W_2}+{\mathcal{D}}_{V_1}
+{\mathcal{D}}_{V_0},\quad\mbox{for }
V_s\in{\mathcal{K}}_{n,s},\ W_s\in{\mathcal{L}}_{n,s}.$$
Another way of stating this is:--
\begin{corollary}\label{gradedstructure}
There is the following canonical isomorphism of vector spaces.
\begin{equation}\label{gradedisomorphism}
{\mathcal{B}}_n\simeq{\mathcal{K}}_{n,0}\oplus{\mathcal{K}}_{n,1}
\oplus\bigoplus_{s=2}^\infty
\left({\mathcal{K}}_{n,s}\oplus{\mathcal{L}}_{n,s}\right).\end{equation}
\end{corollary}

In order the present the algebra structure on ${\mathcal{B}}_n$ we need more
detail on the spaces ${\mathcal{K}}_{n,s}$ and~${\mathcal{L}}_{n,s}$. It is
well-known and given explicitly in~(\ref{confK}), that the space of
conformal Killing vectors on ${\mathbb{R}}^n$ is isomorphic as a Lie algebra to
${\mathfrak{so}}(n+1,1)$. The spaces ${\mathcal{K}}_{n,s}$ and
${\mathcal{L}}_{n,s}$ are irreducible finite-dimensional representations of
${\mathfrak{so}}(n+1,1)$. Specifically,
\begin{equation}\label{Kns}
{\mathcal{K}}_{n,s}\simeq\overbrace{\raisebox{-10pt}{\begin{picture}(100,30)
\put(0,5){\line(1,0){100}}
\put(0,15){\line(1,0){100}}
\put(0,25){\line(1,0){100}}
\put(0,5){\line(0,1){20}}
\put(10,5){\line(0,1){20}}
\put(20,5){\line(0,1){20}}
\put(30,5){\line(0,1){20}}
\put(40,5){\line(0,1){20}}
\put(70,5){\line(0,1){20}}
\put(80,5){\line(0,1){20}}
\put(90,5){\line(0,1){20}}
\put(100,5){\line(0,1){20}}
\put(55,10){\makebox(0,0){$\cdots$}}
\put(55,20){\makebox(0,0){$\cdots$}}
\put(102,7){\makebox(0,0)[l]{$\circ$}}
\end{picture}}}^{\mbox{$s$}}\end{equation}
and 
\begin{equation}\label{Lns}
{\mathcal{L}}_{n,s}\simeq\overbrace{\raisebox{-10pt}{\begin{picture}(100,30)
\put(0,5){\line(1,0){80}}
\put(0,15){\line(1,0){100}}
\put(0,25){\line(1,0){100}}
\put(0,5){\line(0,1){20}}
\put(10,5){\line(0,1){20}}
\put(20,5){\line(0,1){20}}
\put(30,5){\line(0,1){20}}
\put(40,5){\line(0,1){20}}
\put(70,5){\line(0,1){20}}
\put(80,5){\line(0,1){20}}
\put(90,15){\line(0,1){10}}
\put(100,15){\line(0,1){10}}
\put(55,10){\makebox(0,0){$\cdots$}}
\put(55,20){\makebox(0,0){$\cdots$}}
\put(102,17){\makebox(0,0)[l]{$\circ$}}
\end{picture}}}^{\mbox{$s$}}\end{equation}
as Young tableau, where $\circ$ denotes the trace-free part. These isomorphisms
may be derived from results concerning induced modules in representation
theory, namely Lepowsky's generalisation~\cite{lepowsky77} of the
Bernstein-Gelfand-Gelfand resolution. A proof in the language of partial
differential operators appears in~\cite{BCEGprolong}.

{From} now on, let us write ${\mathfrak{g}}$ for the Lie algebra
${\mathfrak{so}}(n+1,1)$. Then
\begin{equation}\label{decompose}
\raisebox{6pt}{${\mathfrak{g}}\otimes{\mathfrak{g}}=\;$}
\raisebox{-25pt}{$\begin{picture}(10,50)(0,-20)
\put(0,5){\line(1,0){10}}
\put(0,15){\line(1,0){10}}
\put(0,25){\line(1,0){10}}
\put(0,5){\line(0,1){20}}
\put(10,5){\line(0,1){20}}
\end{picture}
\raisebox{31pt}{$\;\otimes\;$}
\begin{picture}(10,50)(0,-20)
\put(0,5){\line(1,0){10}}
\put(0,15){\line(1,0){10}}
\put(0,25){\line(1,0){10}}
\put(0,5){\line(0,1){20}}
\put(10,5){\line(0,1){20}}
\end{picture}
\raisebox{31pt}{$\;=\;$}
\begin{picture}(25,50)(0,-20)
\put(0,5){\line(1,0){20}}
\put(0,15){\line(1,0){20}}
\put(0,25){\line(1,0){20}}
\put(0,5){\line(0,1){20}}
\put(10,5){\line(0,1){20}}
\put(20,5){\line(0,1){20}}
\put(22,7){\makebox(0,0)[l]{$\circ$}}
\end{picture}
\raisebox{31pt}{$\;\oplus\;$}
\begin{picture}(27,50)(0,-20)
\put(0,10){\line(1,0){20}}
\put(0,20){\line(1,0){20}}
\put(0,10){\line(0,1){10}}
\put(10,10){\line(0,1){10}}
\put(20,10){\line(0,1){10}}
\put(22,12){\makebox(0,0)[l]{$\circ$}}
\end{picture}
\raisebox{31pt}{$\;\oplus\;{\mathbb R}\;\oplus\;$}
\begin{picture}(20,50)
\put(0,15){\line(1,0){10}}
\put(0,25){\line(1,0){10}}
\put(0,35){\line(1,0){20}}
\put(0,45){\line(1,0){20}}
\put(0,15){\line(0,1){30}}
\put(10,15){\line(0,1){30}}
\put(20,35){\line(0,1){10}}
\put(12,17){\makebox(0,0)[l]{$\circ$}}
\end{picture}
\raisebox{31pt}{$\;\oplus\;$}
\begin{picture}(10,50)(0,-20)
\put(0,5){\line(1,0){10}}
\put(0,15){\line(1,0){10}}
\put(0,25){\line(1,0){10}}
\put(0,5){\line(0,1){20}}
\put(10,5){\line(0,1){20}}
\end{picture}
\raisebox{31pt}{$\;\oplus\;$}
\begin{picture}(10,50)
\put(0,5){\line(1,0){10}}
\put(0,15){\line(1,0){10}}
\put(0,25){\line(1,0){10}}
\put(0,35){\line(1,0){10}}
\put(0,45){\line(1,0){10}}
\put(0,5){\line(0,1){40}}
\put(10,5){\line(0,1){40}}
\end{picture}$}
\end{equation}
and, following~\cite{eastwood05}, we shall write $X\circledcirc Y$ for the
projection 
$${\mathfrak{g}}\otimes{\mathfrak{g}}\ni V\otimes W\longmapsto
V\circledcirc W\in\raisebox{-11pt}{\begin{picture}(25,30)
\put(0,5){\line(1,0){20}}
\put(0,15){\line(1,0){20}}
\put(0,25){\line(1,0){20}}
\put(0,5){\line(0,1){20}}
\put(10,5){\line(0,1){20}}
\put(20,5){\line(0,1){20}}
\put(22,7){\makebox(0,0)[l]{$\circ$}}
\end{picture}}\;.
$$
It is shown in \cite{eastwood05} that the symmetry algebra ${\mathcal{A}}_n$ 
for the Laplacian is isomorphic to the tensor algebra
$\bigotimes{\mathfrak{g}}$ modulo the two-sided ideal generated by 
\begin{equation}\label{joseph}
V\otimes W-V\circledcirc W-\ffrac{1}{2}[V,W]
+\ffrac{n-2}{4n(n+1)}\langle V,W\rangle
\quad\mbox{for }V,W\in{\mathfrak{g}},\end{equation}
where $\langle\enskip,\enskip\rangle$ is the Killing form (normalised in the
usual way, not as in~\cite{eastwood05}). To state the corresponding result for
${\mathcal{B}}_n$ we also need a notation for the projection onto 
\begin{picture}(27,10)(0,10)
\put(0,10){\line(1,0){20}}
\put(0,20){\line(1,0){20}}
\put(0,10){\line(0,1){10}}
\put(10,10){\line(0,1){10}}
\put(20,10){\line(0,1){10}}
\put(22,12){\makebox(0,0)[l]{$\circ$}}
\end{picture}
and we shall write this as $V\otimes W\mapsto V\bullet W$ meaning as an
idempotent homomorphism of ${\mathfrak{g}}\otimes{\mathfrak{g}}$ into itself.
We also need to observe that
$$\begin{picture}(46,10)
\put(0,-2){\line(1,0){40}}
\put(0,8){\line(1,0){40}}
\put(0,-2){\line(0,1){10}}
\put(10,-2){\line(0,1){10}}
\put(20,-2){\line(0,1){10}}
\put(30,-2){\line(0,1){10}}
\put(40,-2){\line(0,1){10}}
\put(42,0){\makebox(0,0)[l]{$\circ$}}
\end{picture}\hookrightarrow
{\mathfrak{g}}\odot{\mathfrak{g}}\odot{\mathfrak{g}}\odot{\mathfrak{g}}
\subset
{\mathfrak{g}}\otimes{\mathfrak{g}}\otimes{\mathfrak{g}}\otimes{\mathfrak{g}}$$
meaning that there is a unique irreducible summand of the symmetric tensor
product ${\mbox{\small$\bigodot^4$}}{\mathfrak{g}}$ of the indicated type. With
these conventions in place, we have:--
\begin{theorem}\label{Algtheo}
The algebra ${\mathcal{B}}_n$ is isomorphic to the tensor algebra
$\bigotimes{\mathfrak{g}}$ modulo the 2-sided ideal generated by 
\begin{equation}\label{gens}
V\otimes W-V\circledcirc W-V\bullet W-\ffrac{1}{2}[V,W]
+\ffrac{(n-4)(n+4)}{4n(n+1)(n+2)}\langle V,W\rangle
\quad\mbox{for }V,W\in{\mathfrak{g}}\end{equation}
and the image of  
$\;\begin{picture}(46,10)
\put(0,-2){\line(1,0){40}}
\put(0,8){\line(1,0){40}}
\put(0,-2){\line(0,1){10}}
\put(10,-2){\line(0,1){10}}
\put(20,-2){\line(0,1){10}}
\put(30,-2){\line(0,1){10}}
\put(40,-2){\line(0,1){10}}
\put(42,0){\makebox(0,0)[l]{$\circ$}}
\end{picture}\,$ in ${\mbox{\small$\bigotimes^4$}}{\mathfrak{g}}$.
\end{theorem}
\noindent As noted in \cite{eastwood05} for ${\mathcal{A}}_n$, we can quotient
firstly by $V\wedge W-\frac12[V,W]$ to realise ${\mathcal{B}}_n$ as a quotient
of ${\mathfrak{U}}({\mathfrak{g}})$, the universal enveloping algebra
of~${\mathfrak{g}}$. Compared to~${\mathcal{A}}_n$, the appearance of
additional generators at $4^{\mathrm{th}}$ order is new. 

In fact, there is a more precise statement from which Theorem~\ref{Algtheo}
easily follows. It appears as Theorem~\ref{generalstory} in
Section~\ref{proofofALG}.

\section{The proof of Theorem~\ref{Utheo}}
\begin{lemma}\label{hilf} Suppose $V^{bcd\cdots ef}$ is a conformal Killing
tensor with $s$ indices. If we define $\phi^{cd\cdots ef}$ according to
\begin{equation}\label{defofphi}
\nabla^{(a}V^{bcd\cdots ef)}=g^{(ab}\phi^{cd\cdots ef)}, 
\end{equation}
then 
\begin{equation}\label{LapV}
\Delta V^{bcd\cdots ef}=(s-1)g^{(bc}\nabla_a\phi^{d\cdots ef)a}-
(n+2s-4)\nabla^{(b}\phi^{cd\cdots ef)}\end{equation}
and 
$$\mbox{\rm the trace-free part of }\nabla^{(a}\nabla^b\phi^{cd\cdots ef)}=0.$$
\end{lemma}
\begin{proof}Taking the trace of (\ref{defofphi}) gives
$$\ffrac{2}{s+1}\nabla_bV^{bcd\cdots ef}=\ffrac{2n}{s(s+1)}\phi^{cd\cdots ef}
+\ffrac{4(s-1)}{s(s+1)}\phi^{cd\cdots ef}$$
whence
\begin{equation}\label{thisisphi}
\phi^{cd\cdots ef}=\ffrac{s}{n+2s-2}\nabla_bV^{bcd\cdots ef}.\end{equation}
If we apply $\nabla_a$ to (\ref{defofphi}) we obtain 
$$\ffrac{1}{s+1}\Delta V^{bcd\cdots ef}
+\ffrac{s}{s+1}\nabla^{(b}\nabla_aV^{cd\cdots ef)a}=
\ffrac{2}{s+1}\nabla^{(b}\phi^{cd\cdots ef)}
+\ffrac{s-1}{s+1}g^{(bc}\nabla_a\phi^{d\cdots ef)a}.$$
In combination with (\ref{thisisphi}), this completes the proof
of~(\ref{LapV}). Since $n\geq 3$ and $s\geq 1$, the coefficient $n+2s-4$ is
always non-zero and the final conclusion now follows by
differentiating~(\ref{LapV}). 
\end{proof}

Now, we are in a position to prove Theorem~\ref{Utheo}. Let us write
$${\mathcal{D}}=T^{abcde\cdots f}
\nabla_a\nabla_b\nabla_c\nabla_d\nabla_e\cdots\nabla_f+
\mbox{\rm lower order terms},$$
where $T^{abcde\cdots f}$ is a non-zero symmetric tensor, namely the symbol
of~${\mathcal{D}}$. This tensor splits uniquely as
$$T^{abcde\cdots f}=V^{abcde\cdots f}+g^{(ab}W^{cde\cdots f)}
+g^{(ab}g^{cd}X^{e\cdots f)}$$
where $V^{abcde\cdots f}$ and $W^{cde\cdots f}$ are symmetric trace-free and 
$X^{e\cdots f}$ is symmetric. By subtracting
$$\Delta^2X^{e\cdots f}\nabla_e\cdots\nabla_f=
X^{e\cdots f}\Delta^2\nabla_e\cdots\nabla_f+\mbox{lower order terms}$$
from ${\mathcal{D}}$ we have found a canonically equivalent symmetry of the
form
$${\mathcal{D}}=V^{abcde\cdots f}
\nabla_a\nabla_b\nabla_c\nabla_d\nabla_e\cdots\nabla_f+
W^{cde\cdots f}\Delta \nabla_c\nabla_d\nabla_e\cdots\nabla_f
+\mbox{\rm lower order terms}$$
and we claim that $V^{abcde\cdots f}$ must be a conformal Killing tensor and 
$W^{cde\cdots f}$ must be a generalised conformal Killing tensor. To see this
we simply compute $\Delta^2{\mathcal{D}}$ and for this task it is convenient to
use the formula
$$\begin{array}{rcl}\Delta^2(fg)&=&f\Delta^2g+4(\nabla^af)\Delta\nabla_ag
+2(\Delta f)\Delta g+4(\nabla^a\nabla^b f)\nabla_a\nabla_b g\\ 
&&\quad{}+4(\Delta\nabla^af)\nabla_a g+(\Delta^2 f)g
\end{array}$$
and its evident extension to tensor expressions. If we write 
$$\begin{array}{crl}{\mathcal{D}}&=&
V^{hijkl\cdots m}\nabla_h\nabla_i\nabla_j\nabla_k\nabla_l\cdots\nabla_m+
W^{jkl\cdots m}\Delta \nabla_j\nabla_k\nabla_l\cdots\nabla_m\\
&&\quad{}+Y^{ijkl\cdots m}\nabla_i\nabla_j\nabla_k\nabla_l\cdots\nabla_m+
Z^{jkl\cdots m}\nabla_j\nabla_k\nabla_l\cdots\nabla_m\\
&&\qquad{}+\mbox{\rm lower order terms}
\end{array}$$
where $Y^{ijkl\cdots m}$ and $Z^{jkl\cdots m}$ are symmetric, then 
$$\begin{array}{rcl}\Delta^2{\mathcal{D}}&=&
V^{hij\cdots m}\Delta^2
\nabla_h\nabla_i\nabla_j\cdots\nabla_m+
W^{j\cdots m}\Delta^3\nabla_j\cdots\nabla_m\\
&&\quad{}+4(\nabla^{(a}V^{hij\cdots m)})
\Delta\nabla_a\nabla_h\nabla_i\nabla_j\cdots\nabla_m+
4(\nabla^{(a}W^{j\cdots m)})\Delta^2\nabla_a\nabla_j\cdots\nabla_m\\
&&\qquad{}+Y^{ij\cdots m}\Delta^2
\nabla_i\nabla_j\cdots\nabla_m\\
&&\quad\qquad{}+\mbox{\rm lower order terms}.
\end{array}$$
Moving the Laplacian to the right hand side of each of these terms gives
$$\Delta^2{\mathcal{D}}={\mathcal{P}}\Delta^2+4(\nabla^{(a}V^{hij\cdots m)})
\nabla_a\nabla_h\nabla_i\nabla_j\cdots\nabla_m\Delta
+\mbox{lower order terms}$$
for some differential operator ${\mathcal{P}}$ and for this to be of the form
$\delta\Delta^2$ for some differential operator $\delta$ forces
(\ref{confKill}) to hold, as required.

To find a constraint on $W^{jik\cdots m}$ we should consider lower order
terms:--
$$\begin{array}{rcl}
\Delta^2{\mathcal{D}}&=&{\mathcal{Q}}\Delta^2
+4(\nabla^{(a}V^{hijkl\cdots m)})
\nabla_a\nabla_h\nabla_i\nabla_j\nabla_k\nabla_l\cdots\nabla_m\Delta\\
&&\quad{}+2(\Delta V^{hijkl\cdots m})
\nabla_h\nabla_i\nabla_j\nabla_k\nabla_l\cdots\nabla_m\Delta\\
&&\qquad{}+4(\nabla^{(a}\nabla^bV^{hijkl\cdots m)})
\nabla_a\nabla_b\nabla_h\nabla_i\nabla_j\nabla_k\nabla_l\cdots\nabla_m\\
&&\qquad\quad{}+4(\nabla^{(a}\nabla^bW^{jkl\cdots m)})
\nabla_a\nabla_b\nabla_j\nabla_k\nabla_l\cdots\nabla_m\Delta\\
&&\qquad\qquad{}+4(\nabla^{(a}Y^{ijkl\cdots m)})
\nabla_a\nabla_i\nabla_j\nabla_k\nabla_l\cdots\nabla_m\Delta\\
&&\qquad\qquad\quad{}+\mbox{lower order terms}.
\end{array}$$
If we write $V^{hijkl\cdots m}$ according to (\ref{defofphi}) and substitute
from~(\ref{LapV}), we obtain
$$\begin{array}{rcl}
\Delta^2{\mathcal{D}}&=&{\mathcal{R}}\Delta^2
-2(n+2s-4)(\nabla^{(h}\phi^{ijkl\cdots m)})
\nabla_h\nabla_i\nabla_j\nabla_k\nabla_l\cdots\nabla_m\Delta\\
&&\quad{}+4(\nabla^{(h}\phi^{ijkl\cdots m)})
\nabla_h\nabla_i\nabla_j\nabla_k\nabla_l\cdots\nabla_m\Delta\\
&&\qquad{}+4(\nabla^{(h}\nabla^iW^{jkl\cdots m)})
\nabla_h\nabla_i\nabla_j\nabla_k\nabla_l\cdots\nabla_m\Delta\\
&&\qquad{}\quad+4(\nabla^{(h}Y^{ijkl\cdots m)})
\nabla_h\nabla_i\nabla_j\nabla_k\nabla_l\cdots\nabla_m\Delta\\
&&\qquad\qquad{}+\mbox{lower order terms}
\end{array}$$
from which we deduce that
\begin{equation}\label{deduce}
2(\nabla^{(h}\nabla^iW^{jkl\cdots m)})+2(\nabla^{(h}Y^{ijkl\cdots m)})=
(n+2s-6)(\nabla^{(h}\phi^{ijkl\cdots m)}).\end{equation}
Passing to the next order, we find
$$\begin{array}{rcl}
\Delta^2{\mathcal{D}}&=&{\mathcal{R}}\Delta^2
+4(\Delta\phi^{ajkl\cdots m)})
\nabla_a\nabla_j\nabla_k\nabla_l\cdots\nabla_m\Delta\\
&&\quad{}+4(\Delta\nabla^{(a}W^{jkl\cdots m)})
\nabla_a\nabla_j\nabla_k\nabla_l\cdots\nabla_m\Delta\\
&&\qquad{}+2(\Delta Y^{ijkl\cdots m})
\nabla_i\nabla_j\nabla_k\nabla_l\cdots\nabla_m\Delta\\
&&\qquad\quad{}+4(\nabla^{(a}\nabla^hY^{ijkl\cdots m)})
\nabla_a\nabla_h\nabla_i\nabla_j\nabla_k\nabla_l\cdots\nabla_m\\
&&\qquad\qquad{}+4(\nabla^{(a}Z^{jkl\cdots m)})
\nabla_a\nabla_j\nabla_k\nabla_l\cdots\nabla_m\Delta\\
&&\qquad\qquad{}\quad+\mbox{\rm lower order terms}
\end{array}$$
from which we deduce that 
$$\mbox{the trace-free part of }\nabla^{(a}\nabla^hY^{ijkl\cdots m)}=0.$$
Together with Lemma~\ref{hilf}, it follows from (\ref{deduce}) that
$$\mbox{the trace-free part of }
\nabla^{(a}\nabla^h\nabla^iW^{jkl\cdots m)}=0,$$
as required. This completes the proof of Theorem~\ref{Utheo}.
\hfill$\square$

\medskip
In principle, computations such as these are all that is needed to find
${\mathcal{D}}_V$ and ${\mathcal{D}}_W$ as in Theorem~\ref{Etheo}. For example,
we may arrange that (\ref{deduce}) holds by taking
$$Y^{ijkl\cdots m}=\ffrac{n+2s-6}{2}\phi^{ijkl\cdots m}$$
and from (\ref{thisisphi}) we see that ${\mathcal{D}}_V$ must take the form
\begin{equation}\label{DVf}
{\mathcal{D}}_Vf=V^{ab\cdots c}\nabla_a\nabla_b\cdots\nabla_cf+
\ffrac{s(n+2s-6)}{2(n+2s-2)}(\nabla_aV^{ab\cdots c})\nabla_b\cdots\nabla_cf
+\cdots.\end{equation}
This direct approach, however, is difficult. Fortunately, there is a much
easier way of constructing the operators ${\mathcal{D}}_V$ and
${\mathcal{D}}_W$ and this is done in the next section.

\section{The ambient metric and the proof of Theorem~\ref{Etheo}}
\label{Etheoprf}
The constructions in this section closely follow those of \cite{eastwood05} and
so we shall be brief. Let us consider the Lorentzian quadratic form
$$\tilde{g}_{AB}x^Ax^B\ =\ 2x^0x^\infty+g_{ab}x^ax^b\ =\ 
\raisebox{13pt}{$(x^0,x^a,x^\infty)$}\!
\left(\begin{array}{ccc}0&0&1\\ 0&g_{ab}&0\\ 1&0&0\end{array}\right)\!\! 
\left(\!\begin{array}c x^0\\ x^b\\ x^\infty\end{array}\!\right)$$
on~${\mathbb{R}}^{n+2}$. If we embed
${\mathbb{R}}^n\hookrightarrow{\mathbb{RP}}_{n+1}$ according to 
$$x^a\mapsto\left[\!\begin{array}c 1\\ x^a\\ -x^ax_a/2\end{array}\!\right],$$
then the action of ${\mathrm{SO}}(n+1,1)$ on ${\mathbb{R}}^{n+2}$ preserves
${\mathcal{N}}$ the null cone of $\tilde{g}_{AB}$ and the corresponding
infinitesimal action of ${\mathfrak{g}}= {\mathfrak{so}}(n+1,1)$ on the space
of null directions gives rise to conformal Killing vectors on ${\mathbb{R}}^n$.
Explicitly, if ${\mathfrak{g}}$ is realised as skew tensors $V^{BQ}$ on
${\mathbb{R}}^{n+2}$ in the usual way, then one may check that
\begin{equation}\label{PhiPsi}
V^{BQ}\mapsto V^b=\Phi_BV^{BQ}\Psi^b{}_Q,\end{equation}
where
\begin{equation}\label{defofPhiandPsi}
\Phi_B=(-x^bx_b/2,x_b,1),\enskip \Psi^b{}_Q=(-x^b,\delta^b{}_q,0),\enskip 
\delta^b{}_q=\mbox{Kronecker delta}.\end{equation}
In other words, 
$$V^{BQ}=
\left(\begin{array}{ccc}V^{00}&V^{0q}&V^{0\infty}\\
V^{b0}&V^{bq}&V^{b\infty}\\ 
V^{\infty 0}&V^{\infty q}&V^{\infty\infty}\end{array}\right)
=\left(\begin{array}{ccc}0&r^q&\lambda\\
-r^q&m^{bq}&s^b\\ 
-\lambda&-s^q&0\end{array}\right)$$
corresponds to the conformal Killing vector
\begin{equation}\label{confK}
V^b=-s^b-m^b{}_qx^q+\lambda x^b+r_qx^qx^b-(1/2)x_qx^qr^b.\end{equation}
As in~\cite{eastwood05}, the formula (\ref{PhiPsi}) generalises
\begin{equation}\label{PhiPhiPsiPsi}
V^{BQ\cdots CRDSET}\mapsto V^{b\cdots cde}=\Phi_B\cdots\Phi_C\Phi_D\Phi_E
V^{BQ\cdots CRDSET}\Psi^b{}_Q\cdots\Psi^c{}_R\Psi^d{}_S\Psi^e{}_T.
\end{equation}
to provide an explicit realisation of the isomorphism (\ref{Kns}). Here, 
$V^{BQ\cdots CRDSET}$ is skew in each pair of indices $BQ,\ldots,CR,DS,ET$,
is totally trace-free, and is such that skewing over any three indices gives
zero. Similarly, we may take
\begin{equation}\label{W}
W^{BQ\cdots CRDE}\in\;\raisebox{-10pt}{\begin{picture}(100,30)
\put(0,5){\line(1,0){80}}
\put(0,15){\line(1,0){100}}
\put(0,25){\line(1,0){100}}
\put(0,5){\line(0,1){20}}
\put(10,5){\line(0,1){20}}
\put(20,5){\line(0,1){20}}
\put(30,5){\line(0,1){20}}
\put(40,5){\line(0,1){20}}
\put(70,5){\line(0,1){20}}
\put(80,5){\line(0,1){20}}
\put(90,15){\line(0,1){10}}
\put(100,15){\line(0,1){10}}
\put(55,10){\makebox(0,0){$\cdots$}}
\put(55,20){\makebox(0,0){$\cdots$}}
\put(102,17){\makebox(0,0)[l]{$\circ$}}
\end{picture}}\end{equation}
as totally trace-free, skew in each pair $BQ,\ldots,CR$, symmetric in $DE$, and
such that skewing over any three indices gives zero. Then (\ref{Lns}) is
realised by
\begin{equation}\label{Wrealised}
W^{BQ\cdots CRDE}\mapsto W^{b\cdots c}=\Phi_B\cdots\Phi_C\Phi_D\Phi_E
W^{BQ\cdots CRDE}\Psi^b{}_Q\cdots\Psi^c{}_R.\end{equation}

Following Fefferman and Graham~\cite{fefferman-graham}, we shall use the term
`ambient' to refer to objects defined on~${\mathbb{R}}^{n+2}$. For example,
there is the ambient wave operator 
$$\tilde\Delta=\tilde g^{AB}\frac{\partial^2}{\partial x^A\partial x^B}$$
where $\tilde g^{AB}$ is the inverse of $\tilde g_{AB}$. Let $r=\tilde
g_{AB}x^Ax^B$ so that ${\mathcal{N}}=\{r=0\}$. Suppose that $g$ is an ambient
function homogeneous of degree $w-2$. A simple calculation gives
\begin{equation}\label{Lrg}\tilde\Delta(rg)=r\tilde\Delta g+2(n+2w-2)g.
\end{equation}
In particular, if $w=1-n/2$, then $\tilde\Delta(rg)=r\tilde\Delta g$.
Therefore, if $f$ is homogeneous of degree $1-n/2$, then 
$\tilde\Delta f|_{\mathcal{N}}$ depends only on~$f|_{\mathcal{N}}$ (since $rg$ 
provides the freedom in extending such a function off~${\mathcal{N}}$). This 
defines a differential operator on ${\mathbb R}^n$ and, as detailed
in~\cite{eastwood-graham91}, one may easily verify that it is the Laplacian. 
This construction is due to Dirac~\cite{dirac35} and the main point is that it
is manifestly invariant under the action of ${\mathfrak{so}}(n+1,1)$. We say
that $\Delta$ is conformally invariant acting on conformal densities of weight
$1-n/2$ on ${\mathbb R}^n$. This ambient construction of the Laplacian is a
simple example of the `AdS/CFT correspondence' in physics. A principal feature of 
this correspondence is that calculations are simplified by doing them ambiently
or `in the bulk'. This feature pervades all that follows. 
 
Invariance may also be viewed as follows. Recall that
${\mathfrak{g}}={\mathfrak{so}}(n+1,1)$ is realised as skew tensors~$V^{BQ}$.
Each gives rise to an ambient differential operator 
\begin{equation}\label{ambientDV}
{\mathfrak{D}}_V=V^{BQ}x_B\frac{\partial}{\partial x^Q}
\quad\mbox{where }x_B=x^A\tilde g_{AB}.\end{equation}
It is easily verified that, for $g$ and $f$ of any homogeneity, 
\begin{equation}\label{commutes}
{\mathfrak{D}}_V(rg)=r{\mathfrak{D}}_Vg\quad\mbox{and}\quad
\tilde\Delta{\mathfrak{D}}_Vf={\mathfrak{D}}_V\tilde\Delta f.\end{equation}
The first of these implies that ${\mathfrak{D}}_V$ induces differential
operators on ${\mathbb R}^n$ for densities of any conformal weight: simply
extend the corresponding homogeneous function on ${\mathcal N}$ into
${\mathbb{R}}^{n+2}$, apply ${\mathfrak{D}}_V$, and restrict back
to~${\mathcal{N}}$. In particular, let us denote by ${\mathcal D}_V$ and
$\delta_V$ the differential operators so induced on densities of weight $1-n/2$
and $-1-n/2$, respectively. Bearing in mind the ambient construction of the
Laplacian, it follows immediately from the second equation of (\ref{commutes})
that $\Delta{\mathcal{D}}_V=\delta_V\Delta$. In other words, the infinitesimal
conformal invariance of $\Delta$ gives rise to the symmetries~${\mathcal{D}}_V$
as differential operators. 

The formula (\ref{ambientDV}) generalises to provide further symmetries. It is 
shown in \cite{eastwood05} that the ambient differential operator
\begin{equation}\label{DeeVee}
{\mathfrak{D}}_V=V^{BQ\cdots CRDSET}x_B\cdots x_Cx_Dx_E
\frac{\partial^s}{\partial x^Q\cdots\partial x^R\partial x^S\partial x^T}
\end{equation}
provides a symmetry of the Laplacian for all 
$$V^{BQ\cdots CRDSET}\in\;\raisebox{-10pt}{\begin{picture}(100,30)
\put(0,5){\line(1,0){100}}
\put(0,15){\line(1,0){100}}
\put(0,25){\line(1,0){100}}
\put(0,5){\line(0,1){20}}
\put(10,5){\line(0,1){20}}
\put(20,5){\line(0,1){20}}
\put(30,5){\line(0,1){20}}
\put(40,5){\line(0,1){20}}
\put(70,5){\line(0,1){20}}
\put(80,5){\line(0,1){20}}
\put(90,5){\line(0,1){20}}
\put(100,5){\line(0,1){20}}
\put(55,10){\makebox(0,0){$\cdots$}}
\put(55,20){\makebox(0,0){$\cdots$}}
\put(102,7){\makebox(0,0)[l]{$\circ$}}
\end{picture}}\quad.$$
The proof of Theorem~\ref{Etheo} for conformal Killing tensors is essentially
contained in the following:--
\begin{proposition}\label{one}
The ambient operator {\rm(\ref{DeeVee})} induces a symmetry of $\Delta^2$.
The symbol of this operator is given by~{\rm(\ref{PhiPhiPsiPsi})}.
\end{proposition}
\begin{proof} Firstly, we need to know the ambient description of~$\Delta^2$.
Iterating (\ref{Lrg}) gives
$$\tilde\Delta^2(rg)=r\tilde\Delta^2g+4(n+2w-4)\tilde\Delta g.$$
In particular, if $w=2-n/2$, then $\tilde\Delta^2(rg)=r\tilde\Delta^2g$.
Therefore, if $f$ is homogeneous of degree $2-n/2$, then 
$\tilde\Delta^2 f|_{\mathcal{N}}$ depends only on~$f|_{\mathcal{N}}$ and it is
shown in~\cite{eastwood-graham91} that the resulting differential operator on
${\mathbb{R}}^n$ is~$\Delta^2$. It is it easily verified that (\ref{commutes}) 
holds more generally for ${\mathfrak{D}}_V$ of the form (\ref{DeeVee}). It
follows that  
\begin{equation}\label{morecommuting}
{\mathfrak{D}}_V(rg)=r{\mathfrak{D}}_Vg\quad\mbox{and}\quad
\tilde\Delta^2{\mathfrak{D}}_Vf={\mathfrak{D}}_V\tilde\Delta^2f\end{equation}
for $f$ and $g$ of any homogeneity. Arguing as for the Laplacian shows that the
operator ${\mathcal{D}}_V$ on ${\mathbb{R}}^n$ obtained from ${\mathfrak{D}}_V$
acting on functions homogeneous of degree $w=2-n/2$ is a symmetry
of~$\Delta^2$. There are some details to be verified to make sure that the
symbol of ${\mathcal{D}}_V$ is given by~(\ref{PhiPhiPsiPsi}). However, similar
verifications are done in \cite{eastwood05} and we leave them to the interested
reader.
\end{proof}
The ambient construction of symmetries from tensors $W^{BQ\cdots CRDE}$ as in
(\ref{W}) is less obvious. The following proposition completes the proof of
Theorem~\ref{Etheo}.
\begin{proposition}\label{two}
For any tensor $W^{BQ\cdots CRDE}$ satisfying the symmetries of 
\,{\rm(\ref{W})} the ambient differential operator
$$W^{BQ\cdots CRDE}x_B\cdots x_C
\Big(x_Dx_E\tilde\Delta-2x_D\frac{\partial}{\partial x^E}\Big)
\frac{\partial^{s-2}}{\partial x^Q\cdots\partial x^R}$$
induces a symmetry ${\mathcal{D}}_W$ of $\Delta^2$ of the form 
$${\mathcal{D}}_W=W^{b\dots c}\Delta\nabla_b\cdots\nabla_c
+\mbox{\rm lower order terms}$$
where $W^{b\cdots c}$ is given by~{\rm(\ref{Wrealised})}.
\end{proposition}
It is not too hard to prove this Proposition by direct calculation along the
lines of Proposition~\ref{one}. There is a difference, however, in that the
analogue of the first equation of (\ref{morecommuting}) holds only for $g$ of
homogeneity $-n/2$. Moreover, for an analogue of the second equation, one needs
to use the ambient operator 
$$W^{BQ\cdots CRDE}x_B\cdots x_C
\Big(x_Dx_E\tilde\Delta+6x_D\frac{\partial}{\partial x^E}\Big)
\frac{\partial^{s-2}}{\partial x^Q\cdots\partial x^R}$$
on the right hand side and, even so, it is valid only for homogeneity $2-n/2$.
Of course, these homogeneities are exactly what we need for $\Delta^2$ but
there is a more satisfactory ambient construction giving rise to exactly the
same symmetries, which we shall defer to the following section. The operators
in this more satisfactory construction enjoy the proper generalisation
of~(\ref{morecommuting}), namely~(\ref{bestcommtues}).

\section{The proof of Theorem~\ref{Algtheo}}\label{proofofALG}
We shall prove Theorem~\ref{Algtheo} by a new method, improving
on~\cite{eastwood05}. As a side effect, we shall obtain a straightforward proof
of Proposition~\ref{two}. 

In the previous section, we found in (\ref{ambientDV}) a linear mapping
\begin{equation}\label{basicrepresentation}
{\mathfrak{so}}(n+1,1)={\mathfrak{g}}\ni V\longmapsto {\mathfrak{D}}_V,
\end{equation}
where ${\mathfrak{D}}_V$ is an ambient differential operator acting on
functions homogeneous of degree $w$ for any~$w$. By dint of
(\ref{commutes}), we obtain an induced linear mapping
$${\mathfrak{g}}\ni V\longmapsto{\mathcal{D}}_V,$$
where ${\mathcal{D}}_V$ is a differential operator acting on conformal
densities of weight~$w$. In fact, it is shown in \cite{eastwood05} that 
\begin{equation}\label{DVonRn}
{\mathcal{D}}_Vf=V^a\nabla_af-\ffrac{w}{n}(\nabla_aV^a)f\end{equation}
where $V^a$ is the conformal Killing vector associated to $V\in{\mathfrak{g}}$ 
according to~(\ref{PhiPsi}). 

The mapping (\ref{basicrepresentation}) immediately extends to the whole
tensor algebra $\bigotimes{\mathfrak{g}}$ by
\begin{equation}\label{geegee}
{\mathfrak{g}}\otimes{\mathfrak{g}}\otimes\cdots\otimes{\mathfrak{g}}\ni
U\otimes \cdots\otimes V=X\mapsto{\mathfrak{D}}_X\equiv
{\mathfrak{D}}_U\cdots{\mathfrak{D}}_V\end{equation}
and extended by linearity. It follows from (\ref{commutes}) that 
\begin{equation}\label{bestcommtues}
{\mathfrak{D}}_X(rg)=r{\mathfrak{D}}_Xg\quad\mbox{and}\quad
\tilde\Delta{\mathfrak{D}}_X={\mathfrak{D}}_X\tilde\Delta\end{equation}
and hence that there is an induced series of operators ${\mathcal{D}}_X$ acting
on densities of weight $w$ on ${\mathbb{R}}^n$ and providing symmetries of
$\Delta^k$ when $w=k-n/2$. Of course, for simple tensors $X$ these operators
are obtained by composing the basic operators ${\mathcal{D}}_V$ for
$V\in{\mathfrak{g}}$. However, usefully to compute even the basic composition
${\mathcal{D}}_U{\mathcal{D}}_V$ from (\ref{DVonRn}) for $U,V\in{\mathfrak{g}}$
is difficult. Though this is done in~\cite{eastwood05}, the simpler approach
adopted there is to compute ${\mathfrak{D}}_U{\mathfrak{D}}_V$ instead. The
object is to see how this composition breaks up under (\ref{decompose}) but,
for this purpose, the following argument is even more straightforward.

We compute
$${\mathfrak{D}}_U{\mathfrak{D}}_V=
U^{BQ}x_B\frac{\partial}{\partial x^Q}V^{CR}x_C\frac{\partial}{\partial x^R}=
U^{BQ}V^{CR}x_Bx_C\frac{\partial^2}{\partial x^Q\partial x^R}+
U^{BQ}V_Q{}^Rx_B\frac{\partial}{\partial x^R}\,,$$
which extends by linearity to give
\begin{equation}\label{DeeX}{\mathfrak{D}}_X=
X^{BQCR}x_Bx_C\frac{\partial^2}{\partial x^Q\partial x^R}+
X^{BQ}{}_Q{}^Rx_B\frac{\partial}{\partial x^R}\,,\quad\mbox{for }
X^{BQCR}\in{\mathfrak{g}}\otimes{\mathfrak{g}}.\end{equation}
We can simply apply this formula to tensors $X$ from each of the summands on
the right hand side of (\ref{decompose}). All $X^{BQCR}$ in
${\mathfrak{g}}\otimes{\mathfrak{g}}$ are skew in $BQ$ and $CR$ but the various
summands of (\ref{decompose}) are characterised as follows. 
$$\begin{picture}(25,30)(0,2)
\put(0,5){\line(1,0){20}}
\put(0,15){\line(1,0){20}}
\put(0,25){\line(1,0){20}}
\put(0,5){\line(0,1){20}}
\put(10,5){\line(0,1){20}}
\put(20,5){\line(0,1){20}}
\put(22,7){\makebox(0,0)[l]{$\circ$}}
\end{picture}\raisebox{10pt}{$\;\longleftrightarrow
\left\{\begin{array}l
X^{BQCR}+X^{BCRQ}+X^{BRQC}=0\\
X^{BQCR}\mbox{ is totally trace-free.}\end{array}\right.$}$$
Therefore, the second term in (\ref{DeeX}) vanishes and ${\mathfrak{D}}_X$ is
given by~(\ref{DeeVee}), as expected. 

Next we have
\begin{equation}\label{symmtracefree}
\raisebox{-6pt}{$\begin{picture}(27,20)(0,5)
\put(0,10){\line(1,0){20}}
\put(0,20){\line(1,0){20}}
\put(0,10){\line(0,1){10}}
\put(10,10){\line(0,1){10}}
\put(20,10){\line(0,1){10}}
\put(22,12){\makebox(0,0)[l]{$\circ$}}
\end{picture}\raisebox{5pt}{$\;\leftrightarrow X^{BQCR}=
W^{BC}\tilde g^{QR}-W^{QC}\tilde g^{BR}-W^{BR}\tilde g^{QC}
+W^{QR}\tilde g^{BC}$}$}\end{equation}
where $W^{BC}$ is symmetric trace-free. Therefore, 
$${\mathfrak{D}}_X=W^{BC}x_Bx_C\tilde\Delta
-2W^{QC}x_Cx^R\frac{\partial^2}{\partial x^R\partial x^Q}
-nW^{BR}x_B\frac{\partial}{\partial x^R}
+W^{QR}r\frac{\partial^2}{\partial x^Q\partial x^R}$$
and, when acting on functions homogeneous of degree~$w$,
$${\mathfrak{D}}_X=W^{BC}\Big(x_Bx_C\tilde\Delta
-(n+2w-2)x_B\frac{\partial}{\partial x^C}
+r\frac{\partial^2}{\partial x^B\partial x^C}\Big).$$
There are two immediate consequences of this formula. Firstly, when $w=1-n/2$,
the appropriate homogeneity for the Laplacian, we obtain
$${\mathfrak{D}}_X=W^{BC}\Big(x_Bx_C\tilde\Delta
+r\frac{\partial^2}{\partial x^B\partial x^C}\Big)$$
and the induced operator ${\mathcal{D}}_X$ on ${\mathbb{R}}^n$ is clearly of
the form ${\mathcal{P}}\Delta$. Therefore, this summand in the decomposition
(\ref{decompose}) is contained in the annihilator ideal. This is confirmed
by~(\ref{joseph}). On the other hand, when $w=2-n/2$, we obtain
$${\mathfrak{D}}_X=W^{BC}\Big(x_Bx_C\tilde\Delta
-2x_B\frac{\partial}{\partial x^C}
+r\frac{\partial^2}{\partial x^B\partial x^C}\Big)$$
and the induced operator on ${\mathbb{R}}^n$ coincides with the statement of
Proposition~\ref{two} in this case. It is also easy to compute the symbol of
the induced operator on~${\mathbb{R}}^n$ as follows.
\begin{lemma}
$$\tilde g^{QR}\Psi^b{}_Q\Psi^c{}_R=g^{bc}\qquad
\Phi_B\tilde g^{BQ}\Psi^b{}_Q=0\qquad
\Phi_B\Phi_C\tilde g^{BC}=0$$
\end{lemma}
\begin{proof} These are simple computations from (\ref{defofPhiandPsi}).
\end{proof}
{From} this lemma, if $X^{BQCR}$ is of the form given 
in~(\ref{symmtracefree}), then 
$$\Phi_B\Phi_CX^{BQCR}\Psi^b{}_Q\Psi^c{}_R=\Phi_B\Phi_CW^{BC}g^{bc}=Wg^{bc}$$
where $W$ is given by~(\ref{Wrealised}). It follows that the induced operator
on ${\mathbb{R}}^n$ is of the form 
$$Wg^{ab}\nabla_a\nabla_b+\mbox{lower order terms}=
W\Delta+\mbox{lower order terms}.$$
Therefore, we have proved Proposition~\ref{two} for second order operators. Let
us return to analysing the effect of the various summands of
${\mathfrak{g}}\otimes{\mathfrak{g}}$ in~(\ref{DeeX}). 

Next we have
$${\mathbb{R}}\leftrightarrow X^{BQCR}=
V\ffrac{1}{n(n+1)(n+2)}
\big(\tilde g^{QC}\tilde g^{BR}-\tilde g^{BC}\tilde g^{QR}\big)\quad
\mbox{for constant $V$}.$$
The normalisation is arranged so that the Killing form 
$\langle\enskip,\enskip\rangle:
{\mathfrak{g}}\otimes{\mathfrak{g}}\to{\mathbb{R}}$ gives
$$\begin{array}{rcl}
X^{BQCR}\longmapsto-nX^{BQ}{}_{BQ}&=&-nV\ffrac{1}{n(n+1)(n+2)}
\big(\delta^Q{}_B\delta^B{}_Q-\delta^B{}_B\delta^Q{}_Q\big)\\[10pt]
&=&-nV\ffrac{1}{n(n+1)(n+2)}
\big((n+2)-(n+2)^2\big)=V.\end{array}$$
We compute
$${\mathfrak{D}}_X=V\ffrac{1}{n(n+1)(n+2)}\big(
x^Qx^R\frac{\partial^2}{\partial x^R\partial x^Q}+r\tilde\Delta
+(n+1)x^R\frac{\partial}{\partial x^R}\big).$$
Therefore, when acting on functions homogeneous of degree~$w$, 
$${\mathfrak{D}}_X=V\ffrac{1}{n(n+1)(n+2)}\big(r\tilde\Delta+w(n+w)\big).$$
Hence, the corresponding action of ${\mathcal{D}}_X$ on ${\mathbb{R}}^n$ is
$${\mathcal{D}}_Xf=\ffrac{w(n+w)}{n(n+1)(n+2)}Vf$$
for conformal densities of weight~$w$. 
In particular, if $w=1-n/2$ then 
$${\mathcal{D}}_Xf=-\ffrac{(n-2)}{4n(n+1)}Vf.$$
This is exactly as predicted in~(\ref{joseph}). If $w=2-n/2$, however, then 
$${\mathcal{D}}_Xf=-\ffrac{(n-4)(n+4)}{4n(n+1)(n+2)}Vf,$$
in agreement with~(\ref{gens}). 

Next we have 
$$\begin{picture}(20,40)(0,10)
\put(0,15){\line(1,0){10}}
\put(0,25){\line(1,0){10}}
\put(0,35){\line(1,0){20}}
\put(0,45){\line(1,0){20}}
\put(0,15){\line(0,1){30}}
\put(10,15){\line(0,1){30}}
\put(20,35){\line(0,1){10}}
\put(12,17){\makebox(0,0)[l]{$\circ$}}
\end{picture}\raisebox{15pt}{$\;\longleftrightarrow
\left\{\begin{array}l
X^{BQCR}+X^{CRBQ}=0\\
X^{BQCR}\mbox{ is totally trace-free.}\end{array}\right.$}$$
In this case both terms in (\ref{DeeX}) evidently vanish. 

Next we have
$$\begin{picture}(10,30)
\put(0,5){\line(1,0){10}}
\put(0,15){\line(1,0){10}}
\put(0,25){\line(1,0){10}}
\put(0,5){\line(0,1){20}}
\put(10,5){\line(0,1){20}}
\end{picture}\raisebox{10pt}{$\;\;\leftrightarrow X^{BQCR}
=\ffrac{1}{2n}\big(V^{BR}\tilde g^{QC}-V^{QR}\tilde g^{BC}
-V^{BC}\tilde g^{QR}+V^{QC}\tilde g^{BR}\big)$}$$
where $V^{BR}$ is skew. The normalisation is arranged so that the Lie bracket
${\mathfrak{g}}\otimes{\mathfrak{g}}\to{\mathfrak{g}}$ gives 
$$X^{BQCR}\longmapsto X^B{}_Q{}^{QR}-X^R{}_Q{}^{QB}=
\ffrac{1}{2}V^{BR}-\ffrac{1}{2}V^{RB}=V^{BR}.$$
Since $X^{BQCR}=-X^{CRBQ}$, the first term in (\ref{DeeX}) vanishes. Therefore,
$${\mathfrak{D}}_X=X^{BQ}{}_Q{}^Rx_B\frac{\partial}{\partial x^R}=
\ffrac{1}{2}V^{BR}x_B\frac{\partial}{\partial x^R}=
\ffrac{1}{2}{\mathfrak{D}}_V.$$
This accounts for the term $\ffrac{1}{2}{\mathfrak{D}}_{[V,W]}$ in both
(\ref{joseph}) and~(\ref{gens}). 

The final summand corresponds to totally skew tensors $X^{BQCR}$ and for these 
it is clear that ${\mathfrak{D}}_X$ given by (\ref{DeeX}) vanishes. It
accounts for the presence of this summand in the ideals defined by
(\ref{joseph}) or~(\ref{gens}). 

In summary, by considering the effect in (\ref{DeeX}) of tensors from the
various summands in the decomposition (\ref{decompose}) of
${\mathfrak{g}}\otimes{\mathfrak{g}}$, we have verified (\ref{joseph})
and~(\ref{gens}). It is also worthwhile recording what we have shown for a
general conformal weight~$w$. 

\begin{theorem}\label{generalstory}
Suppose that $X^a$ and $Y^a$ are conformal Killing vector fields on
${\mathbb{R}}^n$ corresponding to $X^{AP}$ and~$Y^{AP}$, respectively, in
${\mathfrak{g}}={\mathfrak{so}}(n+1,1)$. Then
$${\mathcal{D}}_X{\mathcal{D}}_Yf={\mathcal{D}}_{X\circledcirc Y}f
+{\mathcal{D}}_{X\bullet Y}f+\ffrac{1}{2}{\mathcal{D}}_{[X,Y]}f
+\ffrac{w(n+w)}{n(n+1)(n+2)}{\langle X,Y\rangle}f$$
on densities of weight~$w$. Here, 
\renewcommand{\labelitemi}{\labelitemiii}\begin{itemize}
\item ${\mathcal{D}}_Xf=X^a\nabla_af-\ffrac{w}{n}(\nabla_aX^a)f$.
\item $(X\circledcirc Y)^{ab}=
\ffrac{1}{2}X^aY^b+\ffrac{1}{2}X^bY^a-\ffrac{1}{n}X^cY_cg^{ab}$ is a conformal
Killing tensor and for a general conformal Killing tensor~$V^{ab}$,
$${\mathcal{D}}_Vf=V^{ab}\nabla_a\nabla_bf
-\ffrac{2(w-1)}{n+2}(\nabla_aV^{ab})\nabla_bf
+\ffrac{w(w-1)}{(n+2)(n+1)}(\nabla_a\nabla_bV^{ab})f.$$
\item $X\bullet Y=\ffrac{1}{n}X^aY_a=W$ satisfies
$\nabla_a\nabla_b\nabla_cW=g_{(ab}\phi_{c)}$ and, for such a field in general,
$${\mathcal{D}}_Wf=W\Delta f-\ffrac{n+2w-2}{2}(\nabla^aW)\nabla_af
+\ffrac{w(n+2w-2)}{2(n+2)}(\Delta W)f.$$
\item $[X,Y]^{a}=X^b\nabla_bY^a-Y^b\nabla_bX^a$ is a conformal Killing field.
\item $\langle X,Y\rangle=(\nabla_bX^a)(\nabla_aY^b)
-\ffrac{n-2}{n^2}(\nabla_aX^a)(\nabla_bY^b)
-\ffrac{2}{n}X^a\nabla_a\nabla_bY^b-\ffrac{2}{n}Y^a\nabla_a\nabla_bX^b$ is
constant. 
\end{itemize}
Within
$({\mathfrak{g}}\otimes{\mathfrak{g}})\oplus{\mathfrak{g}}\oplus{\mathbb{R}}$,
however, these operations are defined as in Section~\ref{DefineandState}. 
\end{theorem}
\begin{proof}Apart from $W=X\bullet Y$ the various formulae have just been
established or are taken from~\cite{eastwood05} (with a minor rearrangement for
$\langle X,Y\rangle$ on~${\mathbb{R}}^n$). To complete the proof, therefore, it
remains to establish the formula on ${\mathbb{R}}^n$ for $X\bullet Y$ and the
formula for ${\mathcal{D}}_W$ in general. One possibility is to compute the
composition ${\mathcal{D}}_X{\mathcal{D}}_Y$ in full and collect terms. Though 
this certainly works, there is a short cut as follows. Certainly,
$$\begin{array}{rcl}{\mathcal{D}}_X{\mathcal{D}}_Y&=&
X^aY^b\nabla_a\nabla_b+\mbox{lower order terms}\\
&=&(\ffrac{1}{2}X^aY^b+\ffrac{1}{2}X^bY^a-\ffrac{1}{n}X^cY_cg^{ab})
\nabla_a\nabla_b+\ffrac{1}{n}X^cY_c\Delta+\mbox{lower order terms}\\
&=&(X\circledcirc Y)^{ab}\nabla_a\nabla_b+\ffrac{1}{n}X^cY_c\Delta
+\mbox{lower order terms}.
\end{array}$$
It follows that $X\bullet Y=\ffrac{1}{n}X^aY_b$, as advertised. Rather than
find the lower order terms by direct computation, we claim that they are forced
by invariance under the conformal action of ${\mathfrak{so}}(n+1,1)$. This is
essentially the argument used in \cite[\S5]{eastwood05} to find explicit
formulae for ${\mathcal{D}}_V$ in case of an arbitrary conformal Killing tensor
$V^{bc\cdots d}$ acting on densities of any weight. In our case, the argument
is as follows. We are looking for a differential operator of the form  
\begin{equation}\label{try}
{\mathcal{D}}_Wf=W\Delta f+\alpha(\nabla^aW)\nabla_af+\beta(\Delta W)f
\end{equation}
and it remains to determine the constants $\alpha$ and $\beta$ in order that
such an operator be conformally invariant under flat-to-flat rescalings of the
standard metric on~${\mathbb{R}}^n$. We shall follow the conventions of
\cite{baston-eastwood90} concerning conformal geometry. If 
$\hat g_{ab}=\Omega^2g_{ab}$ is also a flat metric, then
$$\nabla^a\Upsilon_a=-\ffrac{n-2}{2}\Upsilon^a\Upsilon_a
\quad\mbox{for }\Upsilon_a=(\nabla_a\Omega)/\Omega.$$
Now $W$ has conformal weight $2$ and we are supposing that $f$ has conformal
weight~$w$. It follows that
$$\begin{array}{rcl}\hat\nabla_af&=&\nabla_af+w\Upsilon_af\\
\hat\nabla^aW&=&\nabla^aW+2\Upsilon^aW\\
\hat\Delta f&=&\Delta f
+(n+2w-2)\big(\Upsilon^a\nabla_af+\ffrac{w}{2}\Upsilon^a\Upsilon_af\big)\\
\hat\Delta W&=&\Delta W
+(n+2)\big(\Upsilon^a\nabla_aW+\Upsilon^a\Upsilon_aW\big)
\end{array}$$
whence (\ref{try}) satisfies $\hat{\mathcal{D}}_W={\mathcal{D}}_W$ if and only
if
$$\alpha=-\ffrac{n+2w-2}{2}\quad\mbox{and}\quad
\beta=\ffrac{w(n+2w-2)}{2(n+2)},$$
as required.
\end{proof}
\noindent Notice that (\ref{DVtwo}) and (\ref{DWtwo}) are special cases of
Theorem~\ref{generalstory}. Also notice that ${\mathcal{D}}_W=W\Delta$ when
$w=1/n-2$, which explains the absence of $X\bullet Y$ in the
generators (\ref{joseph}) of the annihilator ideal in this case. 

We shall now complete the proof of Proposition~\ref{two} (and,
hence, of Theorem~\ref{Etheo}). 
\begin{lemma}\label{three}
Suppose that
$$\raisebox{11pt}{$X^{BQ\cdots CRDSET}\in\enskip$}
\begin{picture}(10,30)
\put(0,5){\line(1,0){10}}
\put(0,15){\line(1,0){10}}
\put(0,25){\line(1,0){10}}
\put(0,5){\line(0,1){20}}
\put(10,5){\line(0,1){20}}
\end{picture}
\raisebox{11pt}{$\;\odot\cdots\odot\;$}
\begin{picture}(10,30)
\put(0,5){\line(1,0){10}}
\put(0,15){\line(1,0){10}}
\put(0,25){\line(1,0){10}}
\put(0,5){\line(0,1){20}}
\put(10,5){\line(0,1){20}}
\end{picture}
\raisebox{11pt}{$\;\odot\;$}
\begin{picture}(10,30)
\put(0,5){\line(1,0){10}}
\put(0,15){\line(1,0){10}}
\put(0,25){\line(1,0){10}}
\put(0,5){\line(0,1){20}}
\put(10,5){\line(0,1){20}}
\end{picture}
\raisebox{11pt}{$\;\odot\;$}
\begin{picture}(10,30)
\put(0,5){\line(1,0){10}}
\put(0,15){\line(1,0){10}}
\put(0,25){\line(1,0){10}}
\put(0,5){\line(0,1){20}}
\put(10,5){\line(0,1){20}}
\end{picture}
\raisebox{11pt}{$\;=\bigodot^s\!{\mathfrak{g}}
\subset\bigotimes^s\!{\mathfrak{g}}\,.$}$$
In other words, $X^{BQ\cdots CRDSET}$ is skew in each pair of indices
$BQ,\ldots,CR,DS,ET$, has $2s$ indices in total, and is invariant under 
permutations of the paired indices. Then the operator defined by 
{\rm(\ref{geegee})} is 
$$\begin{array}{rcl}
{\mathfrak{D}}_X&=&\displaystyle X^{BQ\cdots CRDSET}x_B\cdots x_Cx_Dx_E
\frac{\partial^s}
{\partial x^Q\cdots\partial x^R\partial x^S\partial x^T}\\[10pt]
&&\quad{}+\ffrac{s(s-1)}{2}X^{BQ\cdots CRDS}{}_S{}^Tx_B\cdots x_Cx_D
\displaystyle 
\frac{\partial^{s-1}}{\partial x^Q\cdots\partial x^R\partial x^T}\\[10pt]
&&\qquad{}+\mbox{\rm lower order terms}.
\end{array}$$
\end{lemma}
\begin{proof} The derivation of (\ref{DeeX}) from (\ref{ambientDV}) is easily
extended by induction.
\end{proof}
\noindent Suppose that $W^{BQ\cdots CRDE}$ satisfies the symmetries of 
(\ref{W}) as in the statement of Proposition~\ref{two}. Generalising 
(\ref{symmtracefree}), let  $X^{BQ\cdots CRDSET}$ be obtained by forming
$$W^{BQ\cdots CRDE}\tilde g^{ST}
-W^{BQ\cdots CRSE}\tilde g^{DT}-W^{BQ\cdots CRSE}\tilde g^{DT}
+W^{BQ\cdots CRST}\tilde g^{DE}$$
and then symmetrising over the paired indices $BQ,\ldots,CR,DS,ET$.
{From} Lemma~\ref{three}, a short calculation gives
$$\begin{array}{rcl}
{\mathfrak{D}}_X&=&\displaystyle W^{BQ\cdots CRDE}x_B\cdots x_C
\Big(x_Dx_E\tilde\Delta-2x_Dx^S\frac{\partial^2}{\partial x^S\partial x^E}
+r\frac{\partial^2}{\partial x^D\partial x^E}\Big)
\frac{\partial^{s-2}}{\partial x^Q\cdots\partial x^R}\\[10pt]
&&\quad{}-(n+2s-4)W^{BQ\cdots CRDE}x^B\cdots x^Cx^D\displaystyle
\frac{\partial^{s-1}}{\partial x^E\partial x^Q\cdots\partial x^R},
\end{array}$$
where $W^{BQ\cdots CRDE}$ being trace-free ensures that there are no lower
order terms. Therefore, when acting on functions homogeneous of degree~$w$, we
find
$${\mathfrak{D}}_X=W^{BQ\cdots CRDE}x_B\cdots x_C
\Big(x_Dx_E\tilde\Delta-(n+2w-2)x_D\frac{\partial}{\partial x^E}
+r\frac{\partial^2}{\partial x^D\partial x^E}\Big)
\frac{\partial^{s-2}}{\partial x^Q\cdots\partial x^R}$$ 
and, in particular, if $w=2-n/2$ then we have completed the proof
of Proposition~\ref{two}. 

It remains to finish the proof of Theorem~\ref{Algtheo}. As
in~\cite{eastwood05}, this is done by considering the corresponding graded
algebra~${\mathrm{gr}}({\mathcal{B}}_n)$. From Corollary~\ref{gradedstructure}
we know the structure of this algebra---as a vector space it is
(\ref{gradedisomorphism}) and its algebra structure arises from its being a
quotient of the tensor algebra $\bigotimes{\mathfrak{g}}$.
Theorem~\ref{generalstory} with $w=2-n/2$ implies that the elements
(\ref{gens}) are contained in the ideal defining ${\mathcal{B}}_n$, namely the
kernel of the mapping $\bigotimes{\mathfrak{g}}\to{\mathcal{B}}_n$. From these
elements alone, the corresponding graded ideal contains
$$V\otimes W-V\circledcirc W-V\bullet W\quad\mbox{for }V,W\in{\mathfrak{g}}.$$
Let us consider the ideal generated by these elements alone, i.e.\ generated by
${\mathcal{I}}_2$ where we have grouped the decomposition (\ref{decompose})
according to 
$$
\begin{picture}(10,30)
\put(0,5){\line(1,0){10}}
\put(0,15){\line(1,0){10}}
\put(0,25){\line(1,0){10}}
\put(0,5){\line(0,1){20}}
\put(10,5){\line(0,1){20}}
\end{picture}
\raisebox{11pt}{$\;\otimes\;$}
\begin{picture}(10,30)
\put(0,5){\line(1,0){10}}
\put(0,15){\line(1,0){10}}
\put(0,25){\line(1,0){10}}
\put(0,5){\line(0,1){20}}
\put(10,5){\line(0,1){20}}
\end{picture}
\raisebox{11pt}{$\;=\;$}
\begin{picture}(25,30)
\put(0,5){\line(1,0){20}}
\put(0,15){\line(1,0){20}}
\put(0,25){\line(1,0){20}}
\put(0,5){\line(0,1){20}}
\put(10,5){\line(0,1){20}}
\put(20,5){\line(0,1){20}}
\put(22,7){\makebox(0,0)[l]{$\circ$}}
\end{picture}
\raisebox{11pt}{$\;\oplus\;$}
\begin{picture}(27,30)
\put(0,10){\line(1,0){20}}
\put(0,20){\line(1,0){20}}
\put(0,10){\line(0,1){10}}
\put(10,10){\line(0,1){10}}
\put(20,10){\line(0,1){10}}
\put(22,12){\makebox(0,0)[l]{$\circ$}}
\end{picture}
\raisebox{11pt}{$\;\oplus\;{\mathcal{I}}_2$.}
$$
In particular, ${\mathcal{I}}_2$ contains ${\mathfrak{g}}\wedge{\mathfrak{g}}$
and so ${\mathrm{gr}}({\mathcal{B}}_n)$ is a quotient of the symmetric tensor
algebra~$\bigodot^s\!{\mathfrak{g}}$. We have just seen how the differential
operators ${\mathcal{D}}_W$ in Proposition~\ref{two} and hence
Theorem~\ref{Etheo} arise---the representation (\ref{W}) is realised as a
specific submodule of $\bigodot^s\!{\mathfrak{g}}$ and, indeed, this is the
unique submodule of this type. Hence, as a vector space, for $s\geq 2$ we may
write 
${\mathcal{K}}_{n,s}\oplus{\mathcal{L}}_{n,s}\subset\bigodot^s\!{\mathfrak{g}}$
and the corresponding symmetry operators are given by the ambient construction
in a uniform fashion (as developed earlier in this section). More specifically,
$$\raisebox{-10pt}{\begin{picture}(100,30)
\put(0,5){\line(1,0){100}}
\put(0,15){\line(1,0){100}}
\put(0,25){\line(1,0){100}}
\put(0,5){\line(0,1){20}}
\put(10,5){\line(0,1){20}}
\put(20,5){\line(0,1){20}}
\put(30,5){\line(0,1){20}}
\put(40,5){\line(0,1){20}}
\put(70,5){\line(0,1){20}}
\put(80,5){\line(0,1){20}}
\put(90,5){\line(0,1){20}}
\put(100,5){\line(0,1){20}}
\put(55,10){\makebox(0,0){$\cdots$}}
\put(55,20){\makebox(0,0){$\cdots$}}
\end{picture}}
\;\subset\textstyle\bigodot^s\!{\mathfrak{g}}$$
consists of those $X^{BQ\cdots CRDSET}$ such that skewing over any three
indices gives zero and then 
$${\mathcal{K}}_{n,s}\oplus{\mathcal{L}}_{n,s}=
\Big\{X\in\raisebox{-10pt}{\begin{picture}(100,30)
\put(0,5){\line(1,0){100}}
\put(0,15){\line(1,0){100}}
\put(0,25){\line(1,0){100}}
\put(0,5){\line(0,1){20}}
\put(10,5){\line(0,1){20}}
\put(20,5){\line(0,1){20}}
\put(30,5){\line(0,1){20}}
\put(40,5){\line(0,1){20}}
\put(70,5){\line(0,1){20}}
\put(80,5){\line(0,1){20}}
\put(90,5){\line(0,1){20}}
\put(100,5){\line(0,1){20}}
\put(55,10){\makebox(0,0){$\cdots$}}
\put(55,20){\makebox(0,0){$\cdots$}}
\end{picture}}
\;\mbox{ s.t.\ trace(trace(}X))=0\Big\}.$$
For convenience, let us write
${\mathcal{K}}_{n,s}\oplus{\mathcal{L}}_{n,s}\equiv{\mathcal{M}}_{n,s}$. {From}
this viewpoint it is easy to see that the two-sided ideal generated by
${\mathcal{I}}_2$ in $\bigotimes{\mathfrak{g}}$ is not big enough to have
(\ref{gradedisomorphism}) as its quotient and the problem is when $s=4$.
Arguing as in \cite{eastwood05}, or more specifically as in
\cite[Theorem~3]{eastwood04}, we would like to show that
\begin{equation}\label{want}{\mathcal{M}}_{n,s}=
\Big({\mathcal{M}}_{n,s-1}\otimes\;
\raisebox{-11pt}{$\begin{picture}(10,30)
\put(0,5){\line(1,0){10}} 
\put(0,15){\line(1,0){10}} 
\put(0,25){\line(1,0){10}}
\put(0,5){\line(0,1){20}} 
\put(10,5){\line(0,1){20}}
\end{picture}$}\;\Big)\;\cap\;
\Big(\;\raisebox{-11pt}{$\begin{picture}(10,30)
\put(0,5){\line(1,0){10}} 
\put(0,15){\line(1,0){10}} 
\put(0,25){\line(1,0){10}}
\put(0,5){\line(0,1){20}} 
\put(10,5){\line(0,1){20}}
\end{picture}$}\;\otimes{\mathcal{M}}_{n,s-1}\Big)
\quad\mbox{for }s\geq 3\end{equation}
but this is not true when $s=4$ and the problem is with traces. More
specifically, it is shown in \cite[Theorem 2]{eastwood04} and, in any case, is
easily verified as in \cite{eastwood05} that
\begin{equation}\label{specialinear}\raisebox{-11pt}{$\begin{picture}(40,30)
\put(0,5){\line(1,0){40}} 
\put(0,15){\line(1,0){40}} 
\put(0,25){\line(1,0){40}}
\put(0,5){\line(0,1){20}} 
\put(10,5){\line(0,1){20}}
\put(20,5){\line(0,1){20}} 
\put(30,5){\line(0,1){20}}
\put(40,5){\line(0,1){20}}
\end{picture}$}\;=
\Big(\;\raisebox{-11pt}{$\begin{picture}(30,30)
\put(0,5){\line(1,0){30}} 
\put(0,15){\line(1,0){30}} 
\put(0,25){\line(1,0){30}}
\put(0,5){\line(0,1){20}} 
\put(10,5){\line(0,1){20}}
\put(20,5){\line(0,1){20}} 
\put(30,5){\line(0,1){20}}
\end{picture}$}\;\otimes\;
\raisebox{-11pt}{$\begin{picture}(10,30)
\put(0,5){\line(1,0){10}} 
\put(0,15){\line(1,0){10}} 
\put(0,25){\line(1,0){10}}
\put(0,5){\line(0,1){20}} 
\put(10,5){\line(0,1){20}}
\end{picture}$}\;\Big)\;\cap\;
\Big(\;\raisebox{-11pt}{$\begin{picture}(10,30)
\put(0,5){\line(1,0){10}} 
\put(0,15){\line(1,0){10}} 
\put(0,25){\line(1,0){10}}
\put(0,5){\line(0,1){20}} 
\put(10,5){\line(0,1){20}}
\end{picture}$}\;\otimes\;
\raisebox{-11pt}{$\begin{picture}(30,30)
\put(0,5){\line(1,0){30}} 
\put(0,15){\line(1,0){30}} 
\put(0,25){\line(1,0){30}}
\put(0,5){\line(0,1){20}} 
\put(10,5){\line(0,1){20}}
\put(20,5){\line(0,1){20}} 
\put(30,5){\line(0,1){20}}
\end{picture}$}\;\Big).\end{equation}
Therefore, we are asking whether a tensor
$X^{BQCRDSET}$ enjoying the symmetries of the left hand side of
(\ref{specialinear}) and such that 
\begin{equation}\label{given}
X^{BQ}{}_{BQ}{}^{DSET}=0\quad\mbox{and}\quad X^{BQCRDS}{}_{DS}=0
\end{equation}
has the property that all its second traces vanish. This is not the case.
Indeed, a counterexample may be constructed from any trace-free symmetric
tensor~$Z^{BCDE}$. Let $\tilde g^{QRST}\equiv \tilde g^{(QR}\tilde g^{ST)}$ and
then
\begin{equation}\label{constructed}
X^{BQCRDSET}=\mbox{skew}(Z^{BCDE}\tilde g^{QRST})\end{equation}
where `skew' means to take the skew part in the index pairs $BQ,CR,DS,ET$
(thus generalising (\ref{symmtracefree}) to $\;\begin{picture}(46,10)
\put(0,-2){\line(1,0){40}}
\put(0,8){\line(1,0){40}}
\put(0,-2){\line(0,1){10}}
\put(10,-2){\line(0,1){10}}
\put(20,-2){\line(0,1){10}}
\put(30,-2){\line(0,1){10}}
\put(40,-2){\line(0,1){10}}
\put(42,0){\makebox(0,0)[l]{$\circ$}}
\end{picture}\,$). It is readily verified that (\ref{given}) are satisfied but 
that 
$$X^{BQC}{}_Q{}^{DSE}{}_S=\mbox{a non-zero multiple of }Z^{BCDE}.$$
The proof of Theorem~\ref{Algtheo} now reduces to the following two facts. The
first is that the tensor $X\in\bigodot^4\!{\mathfrak{g}}$ constructed in
(\ref{constructed}) induces a non-zero multiple of the differential operator
$$W^{BCDE}x_Bx_Cx_Dx_E\Delta^2$$ 
on ${\mathbb{R}}^n$, which we decreed to be equivalent to zero in
Definition~\ref{equivalence}. On the one hand this shows that 
$\;\begin{picture}(46,10)
\put(0,-2){\line(1,0){40}}
\put(0,8){\line(1,0){40}}
\put(0,-2){\line(0,1){10}}
\put(10,-2){\line(0,1){10}}
\put(20,-2){\line(0,1){10}}
\put(30,-2){\line(0,1){10}}
\put(40,-2){\line(0,1){10}}
\put(42,0){\makebox(0,0)[l]{$\circ$}}
\end{picture}\,$ should be included in the annihilator ideal for
${\mathcal{B}}_n$ as stated in Theorem~\ref{Algtheo}. On the other hand,
the second easy fact is that (\ref{want}) is true for $s\not=4$ and this
implies that no further additions to the ideal are necessary. The first fact is
an elementary calculation. Both will be left to the reader.

\end{document}